\documentclass[11pt]{amsart}

\usepackage{epsfig,amsmath}

\newtheorem{thm}{Theorem}[section]
\newtheorem{lem}[thm]{Lemma}
\newtheorem{cor}[thm]{Corollary}
\newtheorem{prop}[thm]{Proposition}
\newtheorem{question}[thm]{Problem}

\theoremstyle{definition}

\newtheorem{note}[thm]{Note}

\newcommand{\R}{\mathbf{R}}
\newcommand{\C}{\mathbf{C}}
\newcommand{\Z}{\mathbf{Z}}

\renewcommand{\S}{\mathbf{S}}
\renewcommand{\tilde}{\widetilde}

\numberwithin{equation}{section}

\begin{document}

\vspace*{-0.5in}
{\small\noindent 
To Appear in \emph{Comment.
Math. Helv.}}
\vspace*{.5in}

\parskip 5pt
\baselineskip 13pt

\def\E{\mathcal{E}}
\def\<{\langle}
\def\>{\rangle}
\def\Sig{\Sigma}
\def\tSig{\widetilde\Sigma}
\def\tX{\widetilde X}
\def\bx{{\bf x}}
\def\bk{{\bf k}}
\def\s{\sigma}
\def\sdot{\dot\sigma}
\def\tq{\tau_Q}
\def\tg{\tilde\gamma}
\def\tG{\tilde\Gamma}
\def\g{\gamma}
\def\eps{\epsilon}

\title[Skew Loops]{Skew Loops and Quadric Surfaces}

\parskip 0pt

\author{Mohammad Ghomi}
\address{Department of Mathematics, University of South Carolina,
Columbia, SC 29208}
\email{ghomi@math.sc.edu}
\urladdr{www.math.sc.edu/$\sim$ghomi}

\author{Bruce Solomon}
\address{Department of Mathematics, Indiana University,
Bloomington, IN 47405}
\email{solomon@indiana.edu}
\urladdr{php.indiana.edu/$\sim$solomon}

\date{First draft May 2000. Last Typeset \today.}
\thanks{The first author was partially supported by the NSF grant
DMS-0204190}.

\begin{abstract}
A skew loop is a closed curve without parallel tangent lines.
We prove: \emph{The only complete surfaces in $\R^3$ with a
point of positive curvature and no
skew loops
are the quadrics.} In particular: \emph{Ellipsoids
are the only closed surfaces without skew loops.} Our efforts
also yield results about skew loops on cylinders and
positively curved surfaces.
\end{abstract}

\maketitle

\section{Introduction}
Here we study the relationship between surfaces in $\R^3$ and
closed curves without parallel tangent lines. Examples of such
curves, which we call \emph{skew loops}, were first
constructed by B. Segre in 1968
\cite{segre:tangents}\footnote{Porter gave an apparently
independent construction in 1970 \cite{porter:note}.} to
disprove a conjecture of H. Steinhaus. Quite recently, Wu
constructed skew loops in every knot class \cite{wu:knots},
and the first author has written down explicit examples on
convex surfaces \cite{ghomi:shadow}\footnote{In
\cite{ghomi:shadow} skew loops were used to solve the
``shadow problem" formulated by H. Wente, which is related to
the stability of constant mean curvature surfaces}. Despite
this general failure of Steinhaus' conjecture, however, Segre
noted that it does hold for loops that lie on ellipsoids,
paraboloids, and certain symmetric cylinders. Here we add
convex hyperboloids to Segre's list, show that certain {\it
a}symmetric cylinders do admit skew loops, and use these
facts to prove that the positively curved quadrics are
actually \emph{characterized} by the absence of skew loops:

\begin{thm}
\label{thm:main}
Let $M$ be a connected $2$-manifold, and $F\colon M\to\R^3$
be a $C^2$ immersion. Suppose that $F$ has positive Gauss
curvature at a point of $M$. Then the following are
equivalent:
\begin{enumerate}
\item $F(M)$ lies on a quadric surface.

\item $F(M)$ contains no $C^2$ skew loops.
\end{enumerate}
In particular, if  $F$ is a \emph{complete} immersion and
admits no $C^2$ skew loops, it is an embedding, and $M$ is
simply connected.
\end{thm}

Any loop on a right cylinder over an open planar curve has a
pair of vertical tangent lines, and hence cannot be skew. So
for purposes of the implication $2\Rightarrow 1$ in Theorem
\ref{thm:main} the assumption of positive curvature at one
point is not superfluous. Moreover, since closed surfaces
(compact 2-manifolds without boundary) always have such a
point, Theorem \ref{thm:main} yields:

\begin{cor}
Ellipsoids are the only closed $C^2$ surfaces immersed in
$\R^3$ which admit {\rm no} $C^2$ skew loops.\qed
\end{cor}

Characterizations of ellipsoids have a long and
rich history  \cite[p. 151]{bonnesen&fenchel:book},
\cite{petty:ellipsoids}, \cite{heil&martini:special}.
Most such theorems, however, are stated and proved within the
class of convex bodies, where the surfaces are a priori
embedded, and topologically spherical. Ours avoids both these
restrictions.

We prove Theorem \ref{thm:main} by developing a sequence of
intermediate results: In Section 2 we use regular homotopy
to show that positively curved surfaces admit no skew
figure-eights (Proposition~\ref{prop:fig8s}). Applying this
fact in Section 3, we then prove that convex quadrics have no
skew loops. This involves a Lorentzian generalization,
following \cite{fenchel:tantrix} and
\cite{solomon:tantrices}, of Jacobi's Theorem on indicatrices
that bisect the sphere \cite[p.407]{spivak:v3}. In Section
\ref{sec:asymmetry} we prove our asymmetric ``cylinder
lemma'' (Proposition \ref{prop:cylinder}): \emph{Any cylinder
with a strictly convex asymmetric base contains a skew loop.}
We then exploit this fact in Section \ref{sec:ghomi}, using a
stretching argument, to show that surfaces without skew loops
have symmetric local cross sections. By a result of W.
Blaschke, this property characterizes quadrics, and thus
gives Theorem \ref{thm:main}.

We conclude with three appendices. The first proves  a result first
stated by Segre, which gives a strong
converse to the asymmetric cylinder lemma mentioned above, but
still leaves the existence of skew loops on certain cylinders
undetermined. We discuss this and other open problems in
Appendix \ref{appendix:problems}, then conclude with a few
historical notes in Appendix \ref{appendix:history}

\section{Preliminaries: Skew Loops and Their Tantrices}
\label{sec:tantrices}
 A $C^{k}$ immersed \emph{loop} is a $C^k$  mapping
$\gamma\colon\S^1\simeq\R/2\pi\to\R^3$ with nowhere-vanishing
velocity $\gamma'$. We say $\gamma$ is \emph{skew} iff
$\,k\ge 1\,$ and
\begin{equation}
\gamma'(t)\times \gamma'(s)\neq 0 \label{eqn:npl}
\end{equation}
for all
distinct $t,\,s\in\R/2\pi$. The
\emph{tantrix} of
$\gamma$ is  the mapping $\tau\colon \S^1\to\S^2$ given by
$\tau(t) : = \gamma'(t)/ \|\gamma'(t)\|.$

\begin{note}
We will frequently use the
following observations: (i)  affine bijections of $\R^3$ map
skew loops to skew loops, and (ii) $\gamma$ is skew iff $\tau(\S^1)$ is
embedded and disjoint from its antipodal reflection, i.e.,
$\tau(t)\neq\pm\tau(s)$ for all
distinct $t,\,s\in\R/2\pi$.
\end{note}

The \emph{curvature} of a $C^2$ immersed loop  is  the speed of
its tantrix ($\|\tau'(t)\|$). In Sections \ref{sec:nonexistence} and
\ref{sec:ghomi}, we need to perturb skew loops while keeping them skew:

\begin{lem}\label{lem:topology}
$C^2$ skew loops with nonvanishing curvature form an open subset in the space of
all $C^2$ immersed loops in $\R^3$, relative to the $C^2$
topology.
\end{lem}

\begin{proof}
Let $\gamma$ be a $C^2$ skew loop with tantrix $\tau$ and nonvanishing curvature.
Then
$\tau$ is 
$C^1$ immersed. Suppose $\tilde\gamma$ is a $C^2$ loop close to $\gamma$ in
the sense of $C^2$ metric on $C^2$ loop space. Then $\tilde\gamma$ has
nonvanishing curvature as well, and therefore has a $C^1$ immersed tantrix
$\tilde\tau$. Further, 
$\tilde
\tau$  is close to
$\tau$ in the $C^1$ metric. So $\tilde \tau$ is embedded, because $\tau$
is embedded, and embeddings are open in $C^1$ immersed loop space
\cite[p. 37]{Hirsch:differential}. Finally, since $\tau$
avoids its antipodal image, it avoids some neighborhood of
that image. So (by the triangle inequality) $\tilde \tau$
avoids its antipodal image as well, and $\tilde\gamma$ is
skew.
\end{proof}

Deformations of loops through immersions---\emph{regular
homotopies}---arise naturally for us since they continuously
deform the tantrix of a loop as well. A basic theorem of H.
Whitney \cite{whitney:homotopy} states that in
$\R^2\simeq\C$, every loop is regularly homotopic to either
the figure-eight
\begin{equation*}
\g_0(e^{it}) := \cos t(1 + {\bf i}\,\sin t) ,
\end{equation*}
or to one of the degree-$k$ circle coverings
 given by
\begin{equation*}
\g_k(e^{it})
:=e^{i\,k\,t}\ ,\quad  k=\pm1,\,\pm2,\,\dots.
\end{equation*}
On $\S^2\simeq\C\cup\{\infty\}$, however, S. Smale
\cite{smale:homotopy} showed that there are just two regular
homotopy classes: that of the figure-eight $\g_0$, and of the
equator $\g_1$. These facts lead to the following lemmas,
useful both here and in Section 3.

\begin{lem}
\label{lem:s2}
Every $C^2$ loop on $\S^2$ is regularly
homotopic in $\S^2$ to its own tantrix.
\end{lem}

\begin{proof}
The  $C^1$ homotopy $h:[0,\pi/2]\times
I\to\S^2$ given by
\[
h(\theta,\,t) := \sigma(t)\,\cos\theta +
\tau(t)\,\sin\theta
\]
deforms any immersed curve $\sigma$ into its tantrix
$\tau:=\sigma'$. To see that $h$ is
regular, recall  the spherical Frenet equation
$\tau' = -\sigma + \kappa_g\,\nu$ where $\kappa_g$
is the geodesic curvature of $\tau$, and  $\nu :=
\sigma\times\tau$.
Setting $\sigma_\theta(t):=h(\theta,t)$, we compute
$$
\sigma'_\theta(t)
\,=\, \sigma'(t)\,\cos\theta + \tau'(t)\,\sin\theta
\,= \,\tau(t)\,\cos\theta + \kappa_g\,\nu(t)\,\sin\theta
-\sigma(t)\,\sin\theta .
$$
Since $\sigma$, $\tau$, and $\nu$ are orthonormal,
$\sigma'_\theta\ne 0$. So $\sigma_\theta$ is an immersion.
\end{proof}

 Lemma~\ref{lem:s2}  implies that
the tantrix of any $C^2$-immersed loop in $\S^2$ is
immersed, a well-known fact \cite{gluck&pan, rosenberg:constant}
that generalizes to loops on any
positively curved surface:

\begin{lem}
\label{obs:immersed}
The tantrix $\tau$ of
\emph{any} $C^2$-immersed curve $\sigma$ on a
positively curved surface $M$ is \emph{immersed} in
$\S^2$.
\end{lem}
\begin{proof}
Parametrize $\sigma$ by arclength, so that $\tau = \sigma'$.
The component of $\tau'$ along a unit normal ${\bf n}$ on $M$
is then given by $ \left(\tau'\right)^\perp =
\left(\sigma''\right)^\perp = \mathsf{k} (\sigma')\,{\bf n}$
where $\mathsf{k}$ denotes normal curvature. Since $M$ is
positively curved, $\mathsf{k}\neq 0$. Hence $\tau'\ne 0$.
\end{proof}

A $C^k$ \emph{figure-eight} $\alpha$ on a surface $M$ is any
$C^k$ loop regularly homotopic to a loop $\beta$ in an open
coordinate disc $\phi:U\to\R^2$, with $\phi\circ\beta=\g_0$
(the ``standard'' figure-eight above). Lemmas \ref{lem:s2}
and \ref{obs:immersed} yield:

\begin{prop}
\label{prop:fig8s}
Let $f:M\to\R^3$ be a $C^2$-immersed, positively curved
surface. Then the tantrix of any figure-eight on $M$ is again
a figure-eight. In particular, $M$ admits no skew
figure-eights.
\end{prop}

\begin{proof}
By definition,  any figure-eight $\alpha\subset M$ is
regularly homotopic to a copy $\beta$ of our ``standard''
figure-eight $\g_0$ in a coordinate disc $U$.  Lemma
\ref{obs:immersed} then implies that the tantrix
$\tau_\alpha$ of $\alpha$ is regularly homotopic to that of
$\beta$: $ \tau_\alpha\sim\tau_\beta\ $. It therefore
suffices to show that $\tau_\beta$ is a figure-eight  on
$\S^2$.

After a regular homotopy of $\beta$ we may assume that $U$ is
so small that $f(U)$ is a graph over one of its tangent
planes. Then, after an affine transformation, $\beta$ lies in
a coordinate disc $U\subset M$ with image $f(U)$
contained in the graph of a convex $C^2$ function
$h_0:D^2\to \R$, where  $D^2\subset\R^2$ is the open unit
disc. We may then realize $\beta$ as a graph, $\,\beta_0$,
over a figure-eight $\g:\S^1\to D^2$:
$$
\beta_0(t) = \g(t) +
h_0\left(\g(t)\right)\bk
$$
where $\bk:=
(0,0,1)$.
We may also assume (dilate  further if necessary)
that the eigenvalues of the hessian $D^2h_0$ lie between 0
and 1 throughout $D$. Now express the southern hemisphere
of $\S^2$ similarly as the graph of a function
$h_1:D^2\to\R$. The eigenvalues of $D^2h_1$
are everywhere at least 1, so the graphs of
the functions
$$
h_\eps(x) := h_0(x) + \eps\,\left(h_1(x)-h_0(x)\right)
$$
give a deformation of $f(U)$ into  $\S^2$ through positively
curved surfaces. By
Lemma \ref{obs:immersed}, the
tantrices of the figure-eights
$
\beta_\eps(t) := \g(t) + h_\eps(\g(t))\,\bk
$
are all immersed. In particular,
$\tau_\beta\sim\tau_{\beta_1}$. By Lemma~\ref{lem:s2},
$\tau_{\beta_1}\sim\beta_1\,$. Thus
$
\tau_\beta \sim\beta_1,
$
which is a figure-eight on $\S^2$.
\end{proof}

\section{Nonexistence of Skew Loops on Quadrics}
\label{sec:nonexistence}

The tantrix of a $C^3$ loop on $\S^2$, if embedded, bisects
the sphere (\cite{fenchel:tantrix},
\cite{solomon:tantrices}). It follows that the tantrix of a
$C^3$ loop on $\S^2$ crosses either itself or its antipodal
image, and hence that $\S^2$ contains no $C^3$ skew loops.
Segre observed that, by affine invariance, this fact extends
to ellipsoids and elliptic paraboloids. Here we sharpen the
argument in \cite{solomon:tantrices} to rule out $C^2$ skew
loops on these same surfaces\footnote{The absence of $C^2$
skew loops on spheres was also established in 1971 by White
\cite{white:immersion}.}, and craft a Lorentzian version that
includes the two-sheeted hyperboloids.
 
\begin{note}
Our methods in this section do not apply to curves that are
only $C^1$. Further, we do not know whether $\S^2$ admits a
skew loop which is $C^1$ but not $C^2$.
\end{note}
\goodbreak

Let $Q$ denote the symmetric
bilinear form on $\R^3$ characterized by
$$
Q(\bx,\bx) = x^2+y^2-z^2,
$$
for all $\bx:=(x,y,z)\in\R^3$. The connected non-singular
level sets of $Q(\bx,\bx)$ are hyperboloids of revolution,
each homothetic to one of the following:
\begin{eqnarray}
\qquad\Sig &:=& \left\{\,\bx\in\R^3\colon Q(\bx,\bx)= -1\,,\
z>0\,\right\}\label{eq:sigma}\quad\text{(hyperboloid of two sheets)}\\
\tSig &:=& \left\{\,\bx\in\R^3\colon Q(\bx,\bx) = +1\,\right\}\quad\text{(hyperboloid of one sheet)}.
\notag
\end{eqnarray}
Differentiating $Q$ along an arc $\s$ immersed in either
$\Sig$ or $\tSig$ gives
\begin{equation}\label{eqn:10}
Q(\s',\s)\equiv 0.
\end{equation}
\goodbreak
\noindent
Thus:
\begin{lem}\label{obs:normal}
Every point $p$ in $\Sig$ or $\tSig$ is \emph{$Q$-normal} to
that surface at $p$.\qed
\end{lem}

Next, parametrize $\Sig$ and $\tSig$ by
$X\colon\R\times(0,\infty)\to \R^3$ and
$\tX\colon\R\times\R\to\R^3$ respectively as follows:
\begin{eqnarray*}
X(u,v) &:=& \big(\cos(u)\,\sinh(v),\ \sin(u)\,\sinh(v),\
\cosh(v)\big) ,\\
\tX(u,v) &:=& \big(\cos(u)\,\cosh(v),\ \sin(u)\,\cosh(v),\
\sinh(v)\big).
\end{eqnarray*}
Since $Q(X_u,X_u),\ Q(X_v,X_v)>0$, and $Q(X_u,X_v)=0$, $Q$
induces a Riemannian metric on $\Sig$\footnote{This is the
well-known hyperbolic metric on $\Sig$.}. So we may define
the \emph{$Q$-tantrix} of an immersed loop $\sigma$ on
$\Sigma$ via
$$
\tq(t) := {\s'(t)\over
\sqrt{Q\left(\s'(t),\,\s'(t)\right)}}.
$$

\begin{note}
\label{rem:Q-tantrix}
Since $Q(\tq,\tq) = +1$, the $Q$-tantrix of a loop on $\Sig$
lies on $\tSig$. Further, $\tq$ is the radial projection of
the (standard) tantrix $\tau$ into $\tSig$. Therefore, much
like $\tau$, the $Q$-tantrix of a skew loop on $\Sig$ is
embedded, and avoids its antipodal image.
\end{note}

In contrast to $\Sig$, $\tSig$ inherits a Lorentzian
structure from $Q$. Indeed, the vectors
\begin{equation}
\label{eqn:frame}
e^+ := {\tX_u\over \cosh(v)}\ ,\qquad
e^- := \tX_v
\end{equation}
form a global frame on $\tSig$, with
\begin{equation}
\label{eqn:30}
Q(e^+,e^+) = +1\ ,\quad
Q(e^-,e^-)=-1\ ,\quad\mbox{and}\quad Q(e^+,e^-) = 0\ .
\end{equation}
If we project out the $Q$-normal direction, the standard
covariant derivative $D$ on $\R^3$ becomes a torsion-free,
$Q$-preserving connection $\nabla$ on $\tSig$. Let $\omega$
denote the corresponding connection 1-form associated to our
frame $\{e^+,\,e^-\}$ by setting
\begin{equation}
\omega(z) := Q(\nabla_z e^+,\,e^-)\ ,\quad\mbox{for all}\ z\in
T\tSig.
\end{equation}
One may verify
that in the local coordinates associated with $\tX$,
\begin{equation}
\label{eqn:sinh}
\omega = -\sinh(v)\,du,
\end{equation}
and that in conjunction with Lemma~\ref{obs:normal},
differentiation of \eqref{eqn:30} yields
\begin{equation}
\label{eqn:50}
\begin{array}{rcccl}
\nabla_ze^+ &=& -Q(\nabla_{z}e^+,\,e^-)\,e^- &=& -\omega(z)\,e^-,\\
\nabla_ze^- &=& +Q(\nabla_{z}e^-,\,e^+)\,e^+ &=& -\omega(z)\,e^+.
\end{array}
\end{equation}

\begin{lem}
\label{lem:noperiod} If a  loop $\alpha$ in $\tSig$ is
the $Q$-tantrix of a $C^2$ loop on $\Sig$, then
$\int_{\alpha} \omega = 0\,$.
\end{lem}

\begin{proof}
Suppose $\alpha =\tau_Q$, the $Q$-tantrix of an arc $\s$
immersed in $\Sig$. Since $\,\tau_Q\,$ is a multiple of
$\s'$,  (\ref{eqn:10}) implies that
$Q(\tq,\s)\equiv 0$. Lemma \ref{obs:normal} then yields that $\s(t)$ is  tangent to
$\tSig$ at $\tq(t)$.  So we may expand $\s$ relative to the frame field given by
(\ref{eqn:frame}). Since $Q(\s,\s)\equiv -1$, and
$\s$ is $C^2$, this uniquely determines a function
$\theta:\S^1\to\R$ such that
\begin{equation*}
\label{eqn:expansion}
\s(t) = \sinh\theta(t)\,e^+ + \cosh\theta(t)\,e^-\ .
\end{equation*}
Note that we evaluate the frame vectors here at $\tq(t)$.
Differentiating the above with respect to $t$,
using (\ref{eqn:50}), yields
\begin{equation*}
\label{eqn:defect} \nabla_{\tq'}\s =
 \big(\theta'
- \omega(\tq')\big)(\cosh\theta\,e^+ + \sinh\theta\,e^-).
\end{equation*}
On the other hand, by  Lemma \ref{obs:normal}, $\tq(t)$ is $Q$-normal to
$\tSig$ at $\tq(t)$. So
\begin{equation*}\label{eqn:parallel}
0=
\big(\sqrt{Q(\s',\s')}\,\tq\big)^\top =
(\s')^\top =
\big(D_{\tq'}\s\big)^\top =
\nabla_{\tq'}\s,
\end{equation*}
which yields that $\omega(\tq')\equiv\theta'$
along $\tq$. But the integral of $\theta'$ vanishes
along $\tq$, since $\theta$ is continuous and
$\tq$ is a loop. Hence
$
\int_{\tq}\omega=0.
$
\end{proof}

We now have the tools we need to prove that \emph{positively
curved quadrics admit no skew loops}, and thereby establish
half of our main theorem.

\begin{proof}[\textbf{Proof of the implication $1 \Rightarrow 2$ of Theorem
\ref{thm:main}}] There are 3 cases:

\emph{Case 1: Hyperboloids.} Each nappe of a hyperboloid of
two sheets is affinely isomorphic to the hyperboloid $\Sig$
defined by \eqref{eq:sigma}. So it suffices to show that
$\Sig$ admits no $C^2$ skew loops. Suppose, toward a
contradiction, that there exists a $C^2$ skew loop
$\s\colon\S^1\to\Sig$, with $Q$-tantrix $\tq$. Since $\Sig$
is diffeomorphic to a plane, and $\sigma$ may not be a
figure-eight (Proposition \ref{prop:fig8s}), Whitney's
theorem forces $\sigma$ to be regularly homotopic to a
$k$-fold tracing $c_k$ of some horizontal circle, $k\ne0$.
The $Q$-tantrix of $c_k$ is then a $k$-fold tracing $\tau_k$
of the circle $z\equiv 0$ in $\tSig$, and since $\Sig$ has
positive curvature, the homotopy $\s\sim c_k$ induces a
regular homotopy $\tau_Q\sim\tau_k$ (Lemma
\ref{obs:immersed}). By Note~\ref{rem:Q-tantrix}, $\tq$ is
embedded, and disjoint from its own antipodal image. The
embeddedness forces $k=1$, and along with the antipodal
disjointness, this means that $-\tq(\S^1)\cup\tq(\S^1)$
bounds an annular domain
$\Omega\subset \tSig$ with $C^1$ boundary.
Combining Stokes' Theorem with
Lemma~\ref{lem:noperiod}, we then get $ \int_{\Omega}\,d\omega
= \int_{\partial\Omega}\,\omega
= 0. $ By  (\ref{eqn:sinh}),  however, $d\omega
= \cosh(v)\,du\,dv$,  a non-vanishing 2-form.
So the integral of $d\omega$
cannot vanish, and we have our contradiction.

\emph{Case 2: Ellipsoids.} All ellipsoids are affinely
equivalent, so we need only check the spherical case, which
has been discussed by Segre \cite{segre:tangents} and White
\cite{white:immersion}. Alternatively, one can proceed as in
Lemma~\ref{lem:noperiod}, replacing $\sinh$ and $\cosh$ by
$\sin$ and $\cos$ respectively. After suitably restricting
their domains, the parametrizations $X$ and $\tX$ for $\Sig$
and $\tSig$ now become patches for $\S^2$. Arguing as in
Lemma~\ref{lem:noperiod}, one then shows that the tantrix of
a loop on $\S^2$ must annihilate the integral of the
corresponding connection form (which is now $\,-\sin v\,du$).
The final argument of Case 1 then goes over almost verbatim,
because $d(\sin v\,du) = \cos v\,du\,dv$ gives the area form
on $\S^2$, except at the poles, which we can avoid with a
slight rotation.

\emph{Case 3: Paraboloids.} By affine equivalence, it
suffices to rule out skew loops on the graph $z=x^2+y^2$. One
easily checks that this paraboloid, call it $P$, can be
$C^2$-approximated arbitrarily well on any compact subset by
an ellipsoid of the form
\begin{equation}
\label{eqn:fr}
x^2 + y^2 +
\left({z\over 2r}-r\right)^2=r^2.
\end{equation}
Further note that, since $P$ has
positive Gaussian curvature,   any
loop on $P$ has nonvanishing curvature. Thus it follows from
Lemma~\ref{lem:topology} that for sufficiently large $r$, any skew loop on $P$ can
be perturbed to form a skew loop on one of the ellipsoids defined by
\ref{eqn:fr} above.  Such a
loop would contradict the result of Case 2, so
$P$ contains no skew loop.
\end{proof}

\section{Asymmetric Convex Cylinders}\label{sec:asymmetry}

When a  $C^k$ loop $\gamma\colon\S^1\to\R^2$ bounds a convex
domain, we say $\Gamma :=\gamma(\S^1)$ is a $C^k$
\emph{oval}.  We say $\Gamma$ is (centrally) \emph{symmetric}
when reflection through a point leaves it invariant.
Otherwise, it is \emph{asymmetric}. We say $\Gamma$ is
\emph{strictly} convex if $\gamma$ is $C^2$ and  its
curvature  never vanishes. Our main aim in this section is to
show:

\begin{prop}[Cylinder Lemma]
\label{prop:cylinder}
The cylinder over any asymmetric, strictly convex $C^2$
oval $\Gamma\subset\R^2$ contains a $C^2$ skew loop with nonvanishing curvature.
\end{prop}

This follows easily once we prove three preliminary results.
Our strategy boils down to the careful analysis of a
classical parametrization: Recall that when $\Gamma$ is
strictly convex, its outward unit normal $n\colon\Gamma\to
\S^1$ is injective. We may therefore employ the \emph{support
parametrization} $\gamma\colon\R\to\R^2$ of $\Gamma$, given by
\begin{equation}\label{eq:gamma}
\gamma(t):=n^{-1}(e^{it}).
\end{equation}
Note that one loses a derivative in passing from $\Gamma$ to
$\g$. When $\Gamma$ is merely $C^2$, this somewhat
complicates the proof that $\gamma$ is an immersion:

\begin{lem}\label{lem:supportfunction}
Let $\Gamma\subset\R^2$ be a strictly convex $C^2$ oval, with
support parametrization $\gamma$. Then $v:=\|\gamma'\|\neq 0$.
Moreover, $\Gamma$ is  symmetric if and only if  $v$  is
$\pi$-periodic.
\end{lem}

\begin{proof}
Define the \emph{support function} of $\Gamma$ via
\begin{equation}\label{eq:supportfunctiondef}
h(t):=\langle e^{it}, \gamma(t) \rangle\quad
\end{equation}
(real inner product). Since $\{e^{it}, ie^{it}\}$ is a basis
for $\R^2$, we then have a $2\pi$-periodic $C^1$ function
$\mu\colon\R\to\R$ such that
\begin{equation}\label{eq:supportfunction:a}
\gamma(t)
= \bigl(h(t) + i\,\mu(t)\bigr)\,e^{it}.
\end{equation}
By \eqref{eq:gamma}, $e^{it}$ is normal to
$\Gamma$ at $\gamma(t)$, so we also have
\begin{equation}\label{eq:supportfunction:b}
\gamma'(t) = v(t)\,i\,e^{it}.
\end{equation}
Now differentiate \eqref{eq:supportfunction:a} and
compare with \eqref{eq:supportfunction:b} to see that
$\mu = h'$ and
\begin{equation}\label{eq:conclusion1}
\gamma(t)=\big(h(t)+i\,h'(t)\big)e^{it}.
\end{equation}
As $\g$ is $\,C^1\,$, this shows that $h$ is  $C^2$.
Further, differentiating \eqref{eq:conclusion1} and using
\eqref{eq:supportfunction:b}, we get
\begin{equation*}\label{conclusion2}
v = h'' + h.
\end{equation*}
We now make indirect use of the curvature formula $\kappa:=
\langle\gamma'',i\gamma'\rangle/ \|\gamma'\|^3$ to show that
$v\neq 0$. If $\gamma$ is $C^2$, one can
differentiate \eqref{eq:supportfunction:b} to evaluate
$\gamma''$, and directly calculate $\kappa = 1/v$. Since
$\Gamma$ is strictly convex, we have $\kappa\neq 0$, and
hence $v\neq 0$, as claimed. Here $\gamma$ is only $C^1$, so
we first approximate $\Gamma$ in $C^2(\S^1,\R^2)$ by a
sequence of $C^3$ ovals $\Gamma_\ell$. The support
parametrization of each $\Gamma_\ell$ will then be $C^2$, so
that for $\Gamma_\ell$, we do have $\kappa_\ell=1/v_\ell$. But
the curvatures $\kappa_\ell$ and speeds $v_\ell$ of the
$\Gamma_\ell$'s converge uniformly to $\kappa$ and $v$
respectively. In the limit, we therefore obtain
$\kappa = 1/v$ as claimed.

To get our final conclusion, suppose that $\Gamma$ is
symmetric about the origin. The reflection $\rho(x)=-x$ then
sends the tangent line at $\gamma(t)$ to some parallel line
tangent to $\Gamma$. Given \eqref{eq:gamma}, the only such
tangency occurs at $\gamma(t+\pi)$. Thus
$$
\gamma(t+\pi) = -\gamma(t)
$$
for all $t\in\R$. By \eqref{eq:supportfunctiondef} this
forces both $h$ and $v = h''+h$ to be $\pi$-periodic.

Conversely, suppose  $v$ is $\pi$-periodic. Then all its odd
Fourier coefficients must vanish. Since $h''+ h = v$, the
same must hold for $h$, modulo a solution $\langle w,
e^{it}\rangle$ of the homogeneous equation $h''+h=0$. By
\eqref{eq:supportfunctiondef}, however, we eliminate this
anomaly if we translate $\Gamma$ by $-w$. Doing so makes $h$
$\pi$-periodic, and by virtue of \eqref{eq:conclusion1},
the oval parametrized by $\g$ is now $\rho$-invariant. The
original (untranslated) oval $\Gamma$ is then symmetric.
\end{proof}
 
We shall define and denote the \emph{even} and \emph{odd}
parts of a function $f\colon\S^1\to\R$ by
$$
f_+(t) := \frac{f(t) + f(t+\pi)}{2},\quad\text{and}\quad
f_-(t) := \frac{f(t) -f(t+\pi)}{2},
$$
respectively, identifying $\S^1$ with $\R/2\pi$ via
$e^{it}\leftrightarrow t$. With this notation, we can give a
simple condition for the skewness of a ``graphical'' loop
$\tg$ on the cylinder over $\Gamma$:

\begin{lem}\label{lem:skew}
Suppose $\Gamma\subset\R^2$ is a strictly convex $C^2$
oval  with support para\-metrization $\gamma$. Let $z\colon
\S^1\to\R$ be $C^1$, and set  $v:=\|\gamma'\|$,
$\mathbf{k}:=(0,0,1)$.  Then
$\tg(t):=\gamma(t)+z(t)\mathbf{k}$ is a skew loop if and only if for
all $t\in\R$, we have
$$
v_+(t)\, z'_+(t)-v_-(t) \,z'_-(t)\neq 0.
$$
Further, if $z$ is $C^2$, then $\tilde\gamma$ has nonvanishing
      curvature.
\end{lem}
\begin{proof}
Expressing $\g'$ as in \eqref{eq:supportfunction:b} above,
and using the identity $i\,e^{i\tau}\times {\mathbf k} =
e^{i\tau}$, we compute that
\begin{eqnarray*}
\tg'(t)\times\tg'(s)
&=& \big(\gamma'(t)\times\gamma'(s)\big)+
\big(z'(s)\gamma'(t)-z'(t)\gamma'(s)\big)\times\mathbf{k}\\
&=& v(t)v(s) e^{it}\times
e^{is}+v(t)z'(s)\,e^{it}-v(s)z'(t)\,e^{is}\\
          &=& v(t) v(s)\sin(t-s)\,\mathbf{k} +v(t) z'(s)\,e^{it}
- v(s) z'(t)\,e^{is}.
\end{eqnarray*}
Note that $\tg$ fails to be skew whenever this quantity
vanishes for some $t$, $s\in \R$, with $t\not\equiv s$ mod
$2\pi$. Since the $\mathbf{k}$ component vanishes only when
$s\equiv t+\pi$ mod $2\pi$, $\tg$ is thus skew if and only if
\begin{equation*}
v(t)z'(t+\pi)+v(t+\pi)z'(t)\neq 0
\end{equation*}
for all $t\in\R$. Now note that for any function
$f\colon\R/2\pi\to\R$, we have the identities
$$
\begin{array}{lcl}
f(t) = f_+(t) + f_-(t)\ ,&f_+(t+\pi) = f_+(t)   \ , \\
f(t+\pi) = f_+(t)-f_-(t)\ ,& f_-(t+\pi) = -f_-(t)\ .
\end{array}
$$
Applying these to $v$ and $z'$ in the preceding formula establishes the first 
conclusion of the lemma.
     Finally, note that, since $\gamma$ is strictly convex,
     $\|\gamma''(t)\|\neq 0$. Thus, if $z$ is $C^2$, then 
 $\|\tilde\gamma''(t)\|\neq 0$ as well. So $\tilde\gamma$ has nonvanishing
curvature.
\end{proof}

The technical result below will provide the key constructive
step in our proof of the cylinder lemma.

\begin{lem}\label{lem:evenandodd}
Let $e$, $o\colon\S^1\to\R$ be continuous functions which are
even and odd respectively, and  suppose that $e+o>0$. Then
either $o\equiv 0$, or we have a continuous function
$\mu\colon\S^1\to\R$, such that:
$$
 ({\rm 1})\;\hbox{$\int_{\S^1}\mu=0$},\qquad
 ({\rm 2})\;\mu \;is\; even,\quad and\quad
({\rm 3})\;e\,\mu>-o^2.
$$
\end{lem}
\begin{proof}
Assume $o\not\equiv 0$, and identify $\S^1$ with $\R/2\pi$ as
usual. To prove the lemma we will construct a continuous
function $\mu\colon[0,\pi]\to\R$ with
$$
\hbox{$({\rm 1}')$ $\int_0^\pi\mu(t)\ dt = 0$},\qquad
\hbox{$({\rm 2}')$ $\mu(\pi)=\mu(0)$},\quad
\hbox{and}\quad
\hbox{$({\rm 3}')$ $e\,\mu > -o^2$ on $[0,\pi]$}\ .
$$
The even extension of this function to all of $\S^1$ then
clearly has the properties (1), (2), and (3) that we seek.

To begin, observe that our hypotheses automatically imply
$\,e>0\,$ throughout $\,\S^1\,$. Otherwise, the evenness of
$e$ would imply $e\le 0$ at both points of some antipodal
pair $t,-t\in\S^1$. Since we assume $e+o>0$ everywhere, this
would force $o>0$ at both $t$ and $-t$, contradicting the
oddness of $o$. We thus have positivity of $e$, which allows
us to define
$$
\tau := \frac{1}{\pi}
\int_0^\pi \left({o(t)^2\over 1+e(t)}\right)\ dt\ >\ 0\ .
$$
Next, note that the zero set of an odd function is both
nonempty, and invariant under reflection through the origin.
After a rotation, we may therefore assume  $o(0)=o(\pi)
= 0$, and define the function we seek:
$$
\mu(t) := \tau - {o^2(t)\over 1 + e(t)}.
$$
Clearly, $\mu$  satisfies $({1}')$. And we arranged
that $o(0)=o(\pi)=0$, so we have $\mu(0) = \mu(\pi) =
\tau$, which gives $({2}')$.
Finally, we obtain $(3')$ by combining
our definition of $\mu$ with the positivity of $e$ and $\tau$:
$$
e(t)\,\mu = e(t)\,\tau - \left({e(t)\over 1 +
e(t)}\right)\,o(t)^2
> - \left({e(t)\over 1 + e(t)}\right)\,o(t)^2
> -o(t)^2.
$$
This proves the Lemma.
\end{proof}

We now prove the main result of this section, our cylinder
lemma:

\begin{proof}[\bf Proof of Proposition \ref{prop:cylinder}]
By Lemma \ref{lem:skew}, it suffices to
 produce a height function
$z\colon\S^1\to\R$ such that, for all
$t\in\R$,
\begin{equation}\label{eq:cylinder:*}
v_+(t) z'_+(t) > v_-(t) z'_-(t),
\end{equation}
where $v$ is the speed of the support parametrization of
$\Gamma$. First, note that our asymmetry hypothesis on
$\Gamma$ combines with Lemma \ref{lem:supportfunction} to
guaranteed that $v$ is not even, and hence $v_-\not\equiv 0$.
Moreover, being odd, $v_-$ has a well-defined antiderivative
on $\S^1$. We form $z_-$ by taking any such antiderivative
and subtracting off its average on $\S^1$. Clearly, this makes
$z_-$ a (non-trivial) odd function, and because
$v_-$ is continuous, $z_-$ is $C^1$.  Since $z'_- = -v_-$,
\eqref{eq:cylinder:*} now becomes
\begin{equation}\label{eq:cylinder:*'}
v_+(t)\,z'_+(t) > - \bigl(v_-(t)\bigr)^2.
\end{equation}
It remains to construct an \emph{even} $C^1$ function
$z_+\colon\S^1\to\R$ whose derivative satisfies
\eqref{eq:cylinder:*'}. Lemma \ref{lem:evenandodd} does
precisely that: Set $e := v_+$, $o:=v_-\not\equiv 0$ there,
and let  $z'_+:=\mu$. Lemma~\ref{lem:supportfunction} ensures
us that $\,e+o = v_++v_- = v>0\,$, so Lemma
\ref{lem:evenandodd} indeed applies. Conclusions (1) and (2)
of the latter now guarantee that $ z'_+$ has an even
antiderivative $z_+$ on $\S^1$, and conclusion (3) reduces to
the key estimate \eqref{eq:cylinder:*'}.
\end{proof}

\section{Quadricity of surfaces without skew loops}
\label{sec:ghomi}

Our first step in this section is to use the existence of
skew loops on asymmetric convex cylinders (Proposition
\ref{prop:cylinder}) to restrict the symmetry of surfaces
without skew loops:

\begin{lem}\label{lem:crossection}
Let $S\subset\R^3$ be a $C^2$ embedded surface without
skew loops. Suppose that there exists a plane $H\subset\R^3$
which meets the interior of $S$ transversely along a strictly
convex oval $\Gamma:=S\cap H$. Then $\Gamma$ is symmetric.
\end{lem}
\begin{proof}
After a rigid motion we may assume that $H$ coincides with
the $xy$-plane. Since $S$ meets $H$ transversely along
$\Gamma$, we may choose $\epsilon>0$ small enough to make
$$
S':=\{\,(x,y,z)\in S: |z|<\epsilon\,\}
$$
a topological annulus transversal to $H$ with
$\partial S'\cap H=\emptyset$.
Let $C$ denote the cylinder perpendicular to $H$ with base
$\Gamma$.  Then $S'$ may be represented as a graph over $C$.
That is, there exists an open neighborhood $A$ of $\Gamma$ in
$C$ and a $C^2$ function $g\colon A\to\R$ such that  $S'=\{ a
+ g(a)\,\nu(a)\,:\ a\in A \}\,$, where  $\nu$ is the outward
unit normal vector field on $S$. Now use the dilatations
$\mu_c\colon\R^3\to\R^3$, defined for each $c\geq 1$ by
$\mu_c(x,y,z) := (x,y,c\,z)$, to define a  1-parameter family
of $C^2$ functions
\[
g_c\colon A\to\R,\qquad\qquad g_c :=
g\circ\mu_{1/c}.
\]
Note that $g_c$ and its derivatives tend to zero uniformly on
$A$ as $c\to\infty$. This follows  from the continuity of $g$
and the chain rule, because $g=0$ on $\Gamma$, while near
$\Gamma$, the derivatives of $g$ are continuous because
$S'$ intersects $H$ transversally.
 
Suppose now that $\Gamma$ is \emph{not}
symmetric. Then Proposition \ref{prop:cylinder}
gives  a $C^2$ skew loop $\gamma\colon\S^1\to C$ with nonvanishing curvature. After
a (shrinking) dilatation, we may assume that
$\gamma(\S^1)\subset A$. For every $c\geq 1$, we may then
define a loop $\gamma_c$ on the affinely stretched surface
$\mu_c(S')$ by setting
\[
\gamma_c(t) := \gamma(t) +
g_c\big(\gamma(t)\big)\,\nu\big(\gamma(t)\big).
\]
Since $g_c\to 0$ uniformly on $\gamma(\S^1)$ along with its
derivatives as $c\to\infty$, we see that $\gamma_c\to\gamma$
in the $C^2$ sense. It then follows, by Lemma
\ref{lem:topology}, that $\gamma_c$ eventually becomes skew.
Thus, for sufficiently large $c>0$, the stretched surface
$\mu_c(S')$ admits a skew loop. As an affine map, however,
$\mu_{c}$ sends skew loops to skew loops. So $S'$ must itself
admit a skew loop---a contradiction.
\end{proof}

By a \emph{convex body} $K\subset\R^3$ we mean a compact
convex subset with nonempty interior. We say planes
$P_1,\,P_2$ are \emph{close} if we can represent them by
linear equations $\langle n_1, x \rangle= h_1$ and $\langle
n_2,x \rangle= h_2$, with $|n_1 - n_2|^2 + |h_1-h_2|^2 <
\epsilon$ for some $\epsilon>0$.

\begin{thm}[Blaschke \cite{blaschke:affine}]
\label{lem:Blaschke}
Let $K\subset\R^3$ be a convex body, whose boundary is $C^2$
near a point $p\in\partial K$. Suppose that whenever
a plane sufficiently close to $T_p{\partial K}$ intersects $K$,
its intersection with $\partial K$ is centrally symmetric.
Then a neighborhood of $p$ in $\partial K$ lies on a quadric
surface.
\end{thm}

Blaschke's result localizes a theorem of Brunn that
characterizes ellipsoids as convex bodies having only
symmetric cross sections \footnote{Olovjanischnikoff (see
\cite{olov:ellipsoid} and \cite[p. 346]{burton:ellipsoid})
proves an even more general version requiring no regularity
at $p$.}. Coupling it to Lemma~\ref{lem:crossection}, we
quickly complete the proof our main theorem.

\begin{proof}[\textbf{Proof of the implication $2\Rightarrow 1$
of Theorem \ref{thm:main}}]
Let $X\subset M$ be the largest open subset whose image $F(X)$
lies on a quadric. Then $X$ is also closed, and $M$ is
connected, so we need only show that $X\neq\emptyset$. To do
so, let $U$ be an open neighborhood of a point $p$ in $M$
where the curvature is positive. We may choose $U$ small
enough so that $S:=F(U)$ is the graph of a function on the
tangent plane $T_{F(p)}\partial K$. Since the curvature is
positive at $p$,  this function has positive definite Hessian
and is therefore convex. So $S$ lies on the boundary of a
convex body $K\subset\R^3$. Since $S$ has positive curvature at
$F(p)$, the tangent plane $T_{F(p)}\partial K$ intersects $K$
only at $F(p)$. This gives an $\epsilon>0$ so that every plane
$H\subset\R^3$ within distance $\epsilon$ of
$T_{F(p)}\partial K$ satisfies $H\cap\partial S=\emptyset$.
Then $\Gamma:=H\cap\partial K$ lies in $S$. Whenever the
intersection is transversal $\Gamma$ is a $C^2$ strictly
convex oval, because $S$ has positive curvature.
Lemma~\ref{lem:crossection} now makes $\Gamma$ symmetric. But
$\Gamma$ was an arbitrary transverse cross-section of $S$ near
$p$, so Blaschke's Theorem (\ref{lem:Blaschke}) forces a
neighborhood of $p$ to lie on a quadric surface.
This completes the proof.
\end{proof}

\appendix
\section{Symmetric cylinders.}\label{sec:appendix}

Proposition \ref{prop:appendix} below gives a strong converse
to the existence of skew loops on asymmetric cylinders
(Proposition \ref{prop:cylinder}). This result was known to Segre \cite{segre:tangents}, but we recount
     a proof for completeness. Let
us agree that an $L$-periodic unit-speed loop $c:\R\to\R^2$
has \emph{arclength symmetry} with respect to a point
$\,p\in\R^2$ if $c(t+L/2) = p-c(t)$ for all $t\in\R$.

\begin{note}
For \emph{embedded} loops, one can show that arclength
symmetry is equivalent to central symmetry. In particular,
Proposition \ref{prop:appendix} holds for cylinders over
embedded centrally symmetric loops. For immersed loops,
arclength symmetry is slightly stronger than central
symmetry, however; centrally symmetric figure-eights, for
instance, admit no arclength-symmetric parametrization.
Indeed, one can put a skew loop on the cylinder over a
centrally symmetric figure-eight. Example: $\tg(t) =
\left(\cos t,\,\sin 2t,\,{t\over \pi}-\left({t\over
\pi}\right)^{15}\right)$. The arclength symmetry condition
below therefore seems essential.
\end{note}

\begin{prop}\label{prop:appendix}
If a $C^1$ loop $\Gamma\subset\R^2$ admits a
parametrization with arclength symmetry, then the cylinder
$S:=\Gamma\times\R\subset\R^3$ admits no skew loops
       which are transverse to the lines in $S$.
\end{prop}

\begin{proof}
Suppose $\Gamma$ has length $L$, and has an $L$-periodic
parametrization $c:\R\to\R^2$ which is arclength-symmetric
about the origin. Let $S:=\Gamma\times\R$, and suppose
$\tg\colon\S^1\to S$ is a $C^1$ loop. 
We may then reparametrize $\tg$ via
\begin{equation}\label{eq:appendix:1}
\tg(t):=c(t)+z(t)\mathbf{k},
\end{equation}
where $z$ is $C^1$ and $nL$-periodic for some $n\in\Z$. By our
symmetry assumption, ${c}(t+L/2)=-{c}(t)$. Hence
\begin{equation}\label{eq:appendix:2}
c'(t+L/2) = -c'(t),
\end{equation}
and
\begin{equation*}
\int_0^{nL}
\big(z'(t)+z'(t+L/2)\big)\,dt
=\big(z(t)+z(t+L/2)\big)\Big|_0^{nL}=0.
\end{equation*}
The mean value theorem for integrals now gives a
$t_0\in[0,nL]$ for which
\begin{equation}\label{eq:appendix:3}
z'(t_0+L/2)=-z'(t_0).
\end{equation}
Equations \eqref{eq:appendix:1},  \eqref{eq:appendix:2}  and
\eqref{eq:appendix:3} therefore combine to yield
$\tg'(t_0+L/2)=-\tg'(t_0)$. This makes the tangent lines to
$\tg$ at $t_0$ and $t_0+L/2$ parallel, a contradiction.
\end{proof}

\section{Open Problems}\label{appendix:problems}
Which surfaces in $\R^3$ admit skew loops? Theorem
\ref{thm:main} settles this question for surfaces with a
point of positive curvature, so it remains to ask:

\begin{question}
Which \emph{nonpositively} curved surfaces admit skew
loops?
\end{question}

Could it be true that the only surfaces without skew loops and
a point of negative curvature are quadric, mirroring Theorem
\ref{thm:main}? If so, it would remain to study flat
surfaces\footnote{Serge Tabachnikov has recently ruled out
skew loops on negatively curved quadrics, and on
simply-connected flat surfaces \cite{tabachnikov:branes}. In
fact, by extending the technique of White
\cite{white:immersion}, he rules out $n$-dimensional
compact skew ``branes'' on all hyperquadrics in $\R^{n+2}$ for
all $n$.}. Complete flat surfaces are generalized cylinders
\cite{hartman&nirenberg}. When embedded and symmetric, these
admit no skew loops, which are transversal to the generators of the
cylinder, by Proposition \ref{prop:appendix}. The main open question
about flat surfaces is then:

\begin{question}
Which \emph{asymmetric} cylinders admit skew loops?
\end{question}

Proposition \ref{prop:cylinder}, shows that strict convexity
is sufficient, and one can show that the tantrix of any
loop on a cylinder whose base has winding number $|\iota|>1$
must self-intersect; these cylinders do not admit skew loops.

Our work raises some regularity questions too. We state one of
them in Note 3.1: \emph{Does $\S^2$ contain a
$C^1$ skew loop}? The regularity of the underlying surface
raises another: Skew loops necessarily have one derivative, so
a version of Theorem \ref{thm:main} in the $C^1$ category
would be fairly optimal with regard to regularity. Our last
question highlights a simple relevant test case:
\begin{question}
Does a cylinder capped by hemispheres admit
skew loops?
\end{question}

\noindent
This surface is $C^1$ and piecewise quadric; if it admits a
skew loop, Theorem \ref{thm:main} is already optimal.

Finally, we remark that when one regards $\R^3$ as
$\mathbf{RP}^3$ minus a plane at infinity, all ellipsoids are
projectively equivalent, not just to each other, but to the
elliptic paraboloids and 2-sheeted hyperboloids too.  The
referee has observed that these are precisely the quadrics on
which we have ruled out skewloops, and hence our results may
extend to $\mathbf{RP^3}$ in an interesting way. We hope to
explore this possibility in a future paper.

\section{Historical Notes}\label{appendix:history}

\small

According to P. Du Val \cite{duval:segre}, H. Steinhaus
conjectured the non-existence of skew loops in 1966  during a
lecture given at Sussex. B. Segre, present at this lecture,
responded by proposing a counterexample in a lecture of his
own the next day. Segre eventually published a corrected
version of his counterexample in 1968 \cite{segre:tangents}.
Porter's version of the construction in 1970
\cite{porter:note} is somewhat more explicit, but Segre's
paper contains many other results, including the
non-existence of skew loops on spheres. To prove the latter
fact, he appeals to a ``\emph{bel teorema}'' published by W.
Fenchel in 1934 \cite{fenchel:tantrix}: \emph{The tantrix of
a spherical curve, if embedded, bisects $\S^2$}. It seems
that this result was absorbed by very few beside Segre. It
immediately implies a well-known theorem of Jacobi on the
normal indicatrix of a space curve, but the many subsequent
references to Jacobi's Theorem we know (e.g.
\cite{spivak:v3}, \cite{chern:maa}, \cite{docarmo:dg}, and
even Fenchel's own 1951 survey \cite{fenchel:bams}!) make no
mention of it. It has since been rediscovered at least twice:
by Avakumovi\'c \cite{avakumovic:tantrix}, and by the second
author \cite{solomon:tantrices}.

The non-existence of skew loops on spheres was also proved by
J.H. White \cite{white:immersion} in 1971 using a
Morse-theoretic argument. Unlike Segre, who notes that the
result extends to ellipsoids and elliptic paraboloids, White
mentions only the sphere. Neither author suggests that
hyperboloids admit no skew loops, nor gives any hint that they
surmised our main result here.

\section*{Acknowledgements}
We thank Ralph Howard, Richard Montgomery, and Serge
Tabachnikov for their interest and helpful comments. We also
thank the University of California at Santa Cruz and Stanford
University for facilitating the early stages of our
collaboration.

\bigskip

\end{document}